\documentclass[12pt]{article}
\usepackage{amssymb,latexsym}
\def\T{\mathfrak{T}}
\def\R{\mathbf{R}}

\def\m{\mathfrak{m}}

\def\C{\mathbf{C}}

\begin{document}
\title{Meromorphic functions with linearly distributed values and
Julia sets of rational functions}
\author{Walter Bergweiler\thanks{Supported by the G.I.F.,
the German--Israeli Foundation for Scientific Research and
Development, Grant G-809-234.6/2003,
the EU Research Training Network CODY,
and the ESF Research Networking Programme HCAA.} and 
Alexandre Eremenko\thanks{Supported by NSF grant DMS-0555279.}}
\maketitle
\begin{abstract} If the preimage of a four-point set
under a meromorphic function belongs to the real line, then
the image of the real line is contained in a circle in the Riemann
sphere. 
We include an application of this result to
holomorphic dynamics: if the Julia set of a rational function is
contained in a smooth curve then it is contained in a circle.
\end{abstract}

If the full preimage of a two-point set under
a rational function belongs to the real line then the function
maps the real line to a circle (on the Riemann sphere). 
We prove an analog of this simple
fact for transcendental functions meromorphic in the
plane.
\vspace{.1in}

\noindent
{\bf Theorem 1.} {\em Let $f$ be a meromorphic function such that the 
preimage of three points belongs to the real line.
Then $f$ maps the real line 
into a circle, unless
\begin{equation}\label{one}
f(z)=L\left(\frac{1-e^{i(c_1z-b_1)}}{1-e^{i(c_2z-b_2)}}\right),
\end{equation}
where $L$ is a fractional-linear transformation, 
and $c_j,b_j$ are real numbers.} 
\vspace{.1in}

\noindent
{\bf Corollary.} {\em 
If the preimage of four points belongs to the real line
then 
$f$ maps the real line into a circle.} 
\vspace{.1in}

The circle of the Corollary does not have to pass
through the four points, as the example
$f(z)=e^{iz}$ shows:
$f$ maps the real line into the unit circle, and the
preimage of $\{0,e^{i\alpha},e^{i\beta},\infty\},\;\alpha,\beta\in\R$,
is on the real line, but in general there is
no circle passing through these four points.

We give an application of Theorem 1 to holomorphic dynamics.
Fatou \cite[no. 56, p. 250]{F} proved the following:

{\em Let $f$ be a rational function, and suppose that some relatively
open subset of the Julia set $J(f)$
is a simple smooth curve.
Then the Julia set
is either a circle or an arc of a circle.}

Here we call a curve smooth
if it has a tangent at every point.
We generalize this result of Fatou:
\vspace{.1in}

\noindent
{\bf Theorem 2.} {\em If a relatively open subset of the Julia
 set of a rational
function is contained in a smooth curve,
then the Julia set is contained in a circle.}
\vspace{.1in}  

This theorem was first proved in \cite{ES} using a different method.
Our paper was inspired by conversations on Theorem~2 with Sebastian
van Strien whom we thank.

Hamilton \cite{H} constructed a rational
function whose Julia set is a Cantor subset of a rectifiable curve,
but does not belong to any circle. 

In the discussion of examples of Julia sets
we assume that the circle in question is the
real line. Our first example is a ``Blaschke product'' 
(a rational function $f$ such that 
both upper and lower halfplanes are invariant under the second iterate
$f^2$).
The Julia sets of such function
can be either a Cantor set on the real line or coincide with the real line. 

The second class of examples consists of certain real polynomials.
Suppose that the Julia set of a polynomial $f$ of degree $d$
is contained in the real line,
and let $[a,b]$ be the convex hull of the Julia set.
Then $f$ is a real function, $f(\{ a,b\})\subset\{ a,b\}$, all
finite critical points of $f$ belong to $[a,b]$, and all critical values
belong to the complement of $(a,b)$. 
It is easy to show that these properties are also sufficient
for the Julia set to lie on the interval $[a,b]$ of the real line.
If all critical values belong to the
set $\{ a,b\}$ then $f$ is conjugate to a Chebyshev polynomial, or to
the negative of
a Chebyshev polynomial. 

There are also examples which are neither
Blaschke products nor polynomials,
the simplest of them is $f(z)=(z^2-c)/(1+\epsilon z),$
where $c<-2$ and $\epsilon$ is real and
small enough. 
It is desirable
to give some classification of rational functions
whose Julia sets are contained in the real line.
\vspace{.1in}

For the proof of Theorem 1 we need the following 
\vspace{.1in}

\noindent{\bf Lemma.} {\em Let $f$ be a meromorphic function
and $a_1,a_2,a_3$ three points in the Riemann sphere
such that all solutions of the equations
$$f(z)=a_j,\quad j=1,2,3$$
are real. Then either $f$ maps the real line into
the circle passing through the $a_j$, or the order of $f$ is at most $1$.}
\vspace{.1in}

\vspace{.1in}

{\em Proof of the Lemma.} We use the standard facts and
notation of Nevanlinna
theory and Tsuji characteristics for a halfplane \cite{GO}.
One could give a more elementary but longer proof by following the method
of Edrei's paper \cite{E}.

Without loss of generality we assume that
$(a_1,a_2,a_3)=(0,1,-1)$. We have to prove that $f$ is either real
or has order at most $1$. Let $g(z)=\overline{f(\overline{z})}$.
As all $a_j$-points of $f$ are real, we conclude that $f$ and $g$ share
$a_j$-points, counting multiplicity.
If $f=g$ then $f$ is real. If $f\neq g$, then we have
$$N(r,0,f-g)\geq N(r,a_1,f)+N(r,a_2,f)+N(r,a_3,f).$$
On the other hand, 
$$N(r,0,f-g)\leq T(r,f-g)+O(1)\leq 2T(r,f)+O(1).$$
Combining these two inequalities, and using the First Main Theorem
of Nevanlinna,  we conclude
that
\begin{equation}\label{m}
T(r,f)\leq m(r,a_1,f)+m(r,a_2,f)+m(r,a_3,f)+O(1).
\end{equation}
As $f$ omits three values in the upper and lower half-planes, we conclude
that its Tsuji characteristics in both half-planes satisfy 
\begin{equation}\label{est}
\m(r,a_j,f)=\T(r,f)+O(1)=O(\log r),
\end{equation}
see \cite[Ch. III, \S3]{GO}.
Now we use the estimate of the Nevanlinna
proximity functions in terms of the Tsuji proximity functions
\cite[Ch. VI, Lemma 5.2]{GO}
$$\int_R^\infty\frac{m(r,a,f)}{r^3}dr\leq\int_R^\infty\frac{\m^*(r,a,f)}{r^2}dr,\quad R>0,$$
where $\m^*$ is the sum of the Tsuji proximity functions of the upper and lower
half-planes. Using (\ref{est}) and (\ref{m}), we conclude that
$$\int_R^\infty\frac{T(r,f)}{r^3}dr=O(R^{-1}\log R),$$
and as $T$ is increasing, we obtain $T(R,f)=O(R\log R)$, that
is the order of $f$ is at most $1$.
This proves the Lemma.
\vspace{.1in}

{\em Proof of Theorem 1.}
Without loss of generality
we assume that the three points are
$0,1,\infty$.
Then we have to prove
that our function $f$ is real, unless it is of the form (\ref{one}).

Put $g(z)=\overline{f(\overline{z})}$.
Then $f$ and $g$ share three values $0,1,\infty$ (counting multiplicity).
We conclude that 
$$f/g=e^u,\quad \mbox{and}\quad (f-1)/(g-1)=e^v,$$
where $u$ and $v$ are some entire functions. 

By solving these equations with respect to $f$ and $g$ we obtain
$$g=\frac{1-e^v}{e^u-e^v}\quad\mbox{and}\quad f=\frac{1-e^v}{1-e^{v-u}}.$$

If $f$ is not real, our Lemma implies that $f$ is of order at
most $1$, so
both functions $u$ and $v$ are affine in this case.
It is easy to see now that any non-constant function
$u$ or $v$ has to be of the form
$i(cz+b)$ where the coefficients $c$ and $b$ are real.

If both $u$ and $v$ are non-constant, or both of them are constant, then
$f$ is of the form (\ref{one}). 

If one of the functions $u,v$ is constant
and the other is not, then $f(z)=L(e^{i(cz+b)})$,
where $L$ is a fractional-linear transformation,
and $c\in\R\backslash\{0\}$,
so $f$ maps the real line into a circle. 

This completes the proof of Theorem 1.
\vspace{.1in}

In the proof of the Corollary
it is enough to consider exceptional functions of the form
$$f(z)=\frac{1-e^{iz}}{1-e^{i(cz+b)}}.$$
where $c\notin\{0,1\}$. (If $c\in\{0,1\}$ then $f(z)=L(e^{iz}),$
where $L$ is a fractional-linear transformation,
and $f$ maps the real line
into a circle). 
Suppose that the preimage of
some $a\in\C\backslash\{0,1,\infty\}$ belongs to
the real line. We are going to prove that $f$ maps the real line into a circle.

So we assume that the equation
\begin{equation}\label{raz}
1-a-e^{iz}+ae^{i(cz+b)}=0
\end{equation}
has only real zeros. The left hand side of this equation is an entire function
of order $1$ with only real zeros, so by Hadamard's factorization 
theorem it is a
product of a real entire function $w$ and
 $e^{-i\gamma z}$ with some real $\gamma$.
Thus the ratio of the left hand side of
(\ref{raz}) and $e^{-i\gamma z}$ is a real
entire function $w$, that is $\overline{w(\overline{z})}=w(z)$ which gives
\begin{eqnarray*}
w(z)&\equiv& 
(1-a)e^{i\gamma z}-e^{i(1+\gamma)z}+ae^{ib}e^{i(c+\gamma)z}\\
&\equiv&(1-\overline{a})e^{-i\gamma z}-e^{-i(1+\gamma)z}+
\overline{a}e^{-ib}e^{-i(c+\gamma)z}.
\end{eqnarray*}
Such an identity can only hold for trivial reasons, because any set of
functions $e^{i\alpha z}$ is linearly independent.
In particular, we conclude that $$\gamma\in\{0,-1-\gamma,-c-\gamma\}.$$
We examine all three possibilities.

1. $\gamma=0$ implies $a=\overline{a},$ then $c=-1$ and $ae^{ib}=-1$. As we assume that
$a\neq 1$ and $b$ is real, the only possibility is that $a=-1$ and $e^{ib}=1.$ Then
$f(z)=(1-e^{iz})/(1-e^{-iz})=-e^{iz}$ maps the real line into the unit circle.

2. $\gamma=-1-\gamma$ implies $\gamma=-1/2$ and then $c\in\{0,1\}$,
the cases
we excluded before.

3. $\gamma=-c-\gamma$ implies $\gamma=-c/2$ and then again $c\in\{0,1\}$,
the cases
excluded before.

This completes the proof of the Corollary.

\vspace{.1in}

{\em Proof of Theorem 2.} Let $V$ be an open set such that the
intersection $J\cap V$ is non-empty and is
contained in a smooth simple curve 
$\gamma$. Let $p\in J\cap V$ be a repelling periodic point.
Replacing $f$ by some iterate we may assume without loss of generality
that $p$ is fixed by $f$.

Consider the Poincar\'e function $F$ associated with $f$ and $p$.
This function satisfies 
$$F(\lambda z)=f(F(z)),\quad F(0)=p,\quad F'(0)=1,$$
and is meromorphic in $\C$.
Here $\lambda=f'(p),\; |\lambda|>1$.

Let $I=F^{-1}(J)$ be the preimage of the Julia set under $F$.
Then $\lambda I=I$. As $F$ is conformal at $0$,
there is a $\lambda^{-1}$-invariant
neighborhood $U$ of $0$ such that $I\cap U$ belongs to a smooth simple
curve $\Gamma$ passing through $0$.
This curve is defined by the property that $F(\Gamma)\subset \gamma$.  
It is easy to see that whenever a subset of a smooth curve is invariant
with respect to multiplication by $\lambda^{-1}$, this subset has
to belong to a straight line through $0$ and $\lambda$ has to be real.
(This argument was used by Fatou \cite[no. 46, p. 229]{F}).
Thus the intersection of $I$ with a neighborhood of $0$
belongs to a straight line. Then $\lambda$-invariance of $I$ implies that
the whole set $I$ is contained in a straight line.
As $F$-preimages of all points of the Julia set
belong to $I$, an application of the Corollary of Theorem~1
shows that the Julia set is contained in a circle.

{\em  

W. B.: Christian-Albrechts-Universität zu Kiel,

bergweiler@math.uni-kiel.de
\vspace{.1in}

A. E.: Purdue University,

eremenko@math.purdue.edu

\begin{thebibliography}{11}
\bibitem{E} A. Edrei, Meromorphic functions with three radially
distributed values, Trans. Amer. Math. Soc., 78 (1955) 276--293.
\bibitem{ES} A. Eremenko and S. van Strien, Rational functions with real
multipliers, manuscript in preparation.
\bibitem{F} P. Fatou, Sur les \'equations fonctionnelles
(Troisi\`eme M\'emoire), Bull. Soc. Math. France, 48 (1920) 208--314.
\bibitem{GO} A. Goldberg and I. Ostrovskii,
Distribution of values of meromorphic functions
(Russian) Moscow, Nauka, 1970.
English translation: Amer. Math. Soc., Providence RI 2008.
\bibitem{H} D. Hamilton, Rectifiable Julia curves, J. London Math. Soc.,
54 (1996) 530--540.
\end{thebibliography}
\end{document}